\documentclass[12pt]{article}
\usepackage{graphicx}
\usepackage{amsmath}
\usepackage{amsfonts}
\usepackage{amsthm}

\newtheorem{theorem}{Theorem}[section]
\newtheorem{lemma}[theorem]{Lemma}
\newtheorem{corollary}[theorem]{Corollary}

\begin{document}

\title{9 Divides no Odd Fibonacci}

\author{Tanya Khovanova}
\date{December 20, 2007}
\maketitle

\begin{abstract}
I discuss numbers that divide no odd Fibonacci. Number 9 plays a special role among such numbers.
\end{abstract}

%%%%%%%%%%%%%%%%%%%%%%%%%%%%%%%%%%%%%%%%%%%%%%%%%%%%%%%%%%%%%%%%%%%%%%%%%%%%%%%%

\section{Introduction}

I stumbled upon the following sentence in the MathWorld
article on the Fibonacci numbers \cite{MWFN}: ``No odd Fibonacci number is
divisible by 17.'' I started wondering if there are other, similar numbers. Of course there are --- no odd Fibonacci number is
divisible by 2. But then, an odd number need not be a Fibonacci number in order
not to be divisible by 2.

So, let us forget about 2 and think about odd numbers. How do we know
that the infinite Fibonacci sequence never produces an odd number that is
divisible by 17? Is 17 the only such odd number? Is 17 the
smallest such odd number? If there are many such odd numbers, how do we
calculate the corresponding sequence?

\section{No odd Fibonacci is divisible by 17}

We will start with a general question: How can we approach puzzles about
the divisibility of Fibonacci numbers? Suppose $K$ is an
integer. Consider the sequence $a_K(n) = F_n \pmod{K}$, of
Fibonacci numbers modulo $K$. The cool thing about this sequence
is that it is periodic. If this is not immediately obvious to you, think
of what happens when a pair of consecutive numbers in the sequence
$a_K(n)$ gets repeated. As a bonus for thinking you will get an
upper bound estimate for this period.

Let us denote the period of $a_K(n)$ by $P_K$. By the way, this
period is called a Pisano
period (see wiki \cite{Pisano}). From the periodicity and the fact that $a_K(0) = 0$, we
see right away that there are infinitely many Fibonacci numbers
divisible by $K$. Are there odd numbers among them? If we trust
MathWorld, then all of the infinitely many Fibonacci numbers divisible
by 17 will be even.

How do we examine the divisibility by $K$ for odd Fibonacci
numbers? Let us look at the Fibonacci sequence modulo 2. As we just
proved, this sequence is periodic. Indeed, every third Fibonacci
number is even. And the evenness of a Fibonacci number is equivalent
to this number having an index divisible by 3.

Now that we know the indices of even Fibonacci numbers we can come
back to the sequence $a_K(n)$. In order to prove that no odd
Fibonacci number is divisible by $K$, it is enough to check that
all the zeroes in the sequence $a_K(n)$ have indices divisible by
3. We already have one zero in this sequence at index 0, which is divisible by 3. Because the sequence $a_K(n)$ is periodic, it will start repeating itself at $a_K(P_K)$. Hence, we need to check that $P_K$ is divisible by 3 and all the zeroes up to $a_K(P_K)$ have indices divisible by 3. When $K =$
17 it is not hard to do the calculations manually.  If you like, try
this exercise. To encourage (or perhaps to discourage) you, here's an
estimate of the scope of the work for this exercise: the Pisano period
for $K =$ 17 is 36.

After I checked that no odd Fibonacci number is ever divisible by 17,
I wanted to find the standard solution for this statement and followed
the trail in MathWorld.  MathWorld sent me on a library trip where I
found the proof of the statement in the book Mathematical
Gems III by Ross Honsberger \cite{Honsberger}.  There was a proof there alright, but it
was tailored to 17 and didn't help me with my questions about other
such odd numbers.

\section{Non-divisibility of odd Fibonacci numbers}

The method we developed for 17 can be used to check any other
number. I trusted this task to my computer. To speed up my program, I
used the fact that the Pisano period for $K$ is never more than
$6K$. Here is the sequence calculated by my trustworthy computer,
which I programmed with, I hope, equal trustworthiness:

\begin{itemize}
\item A133246: Odd numbers n with the property that no odd Fibonacci number is
divisible by n.  \newline
9, 17, 19, 23, 27, 31, 45, 51, 53, 57, 61, 63, 69, 79, 81, 83, 85, 93, 95, 99, $\dots$.
\end{itemize}

The sequence shows that 9 is the smallest odd number that no odd
Fibonacci is ever divisible by, and 17 is the smallest odd
prime with this same property. Here is a trick question for you: Why
is this property of 17 more famous than the same property of 9?

Let us look at the sequence A133246 again. Is this sequence infinite?
Obviously, it should include all multiples of 9 --- hence, it is
infinite. What about prime numbers in this sequence? Is there an
infinite number of primes such that no odd Fibonacci number is
divisible by them? While I do not know the answer, it's worth
investigating this question a little bit further.

\section{Non-divisibility of odd Fibonacci numbers by primes}

From now on, let $K$ be an odd prime. Let us look at the zeroes of
the sequence $a_K(n)$ more closely. Suppose a zero first appears at
the $m$-th place of $a_K(n)$. Then $a_K(m+1) = a_K(m+2) = a$, where $a \ne 0$. In
this case the sequence starting from the $m$-th place is proportional
modulo $K$ to the sequence $a_K(n)$ starting from the 0-th
index. Namely, $a_K(n+m) = a*a_K(n) \pmod{K}$. As $a$ is mutually
prime with $K$, then $a_K(n+m) = 0$ iff $a_K(n) = 0$. From
here, for any index $g$ that is a multiple of $m$, $a_K(g) =
0$. Furthermore, there are no other zeroes in the sequence
$a_K(n)$. Hence, the appearances of 0 in the sequence $a_K(n)$ are
periodic with period $m$.

By the way, $m$ is called a fundamental period; and we just
proved that the Pisano period is a multiple of the fundamental period
for prime $K$.  Hence, the fact that no odd Fibonacci number is
divisible by $K$ is equivalent to the fact that the fundamental
period is not divisible by 3. This is like saying that if the smallest
positive Fibonacci number divisible by an odd prime $K$ is even,
then no odd Fibonacci number is divisible by $K$. In particular, the first Fibonacci number divisible by 17 is $F_9 = 34$; and it is even. Thus we get another proof that 17 divides no odd Fibonacci.

If the remainder of the fundamental period modulo 3 were random, we
would expect that about every third prime number would not 
divide any odd Fibonacci. In reality there are 561 such
primes among the first 1,500 primes (including 2). This is somewhat
more than one third. This gives me hope that there is a non-random
reason for such primes to exist. Consequently, it may be possible to
prove that the sequence of prime numbers that do not divide odd
Fibonacci numbers is infinite.

Can you prove that? Here is the start of this sequence:

\begin{itemize}
\item A133247: Prime numbers p with the property that no odd Fibonacci number is divisible by p. \newline
2, 17, 19, 23, 31, 53, 61, 79, 83, 107, 109, 137, 167, 173, 181, 197, $\dots$.
\end{itemize}

\section{Base sequence}

When I proudly showed to my sons, Alexey Radul and Sergei Bernstein, the two sequences A133246 and A133247 above that I have submitted to the Online Encyclopedia of Integer Sequences (OEIS) \cite{OEIS} they both told me (independently of each other): ``Now you should calculate the base sequence.''

The base sequence means that out of all the numbers that no odd Fibonacci divides we remove multiples of other such numbers. In particular, the base sequence contains all the primes. When I calculated this sequence the result was the following:

\begin{itemize}
\item 2, 9, 17, 19, 23, 31, 53, 61, 79, 83, 107, 109, 137, 167, 173, 181, 197, $\dots$.
\end{itemize}

If you compare the base sequence with the prime sequence A133247, you will see that the only difference is that number 9 belongs to the base sequence.

\begin{theorem}
Number 9 is the only composite in the base sequence.
\end{theorem}

\begin{proof}
Let us denote the fundamental period corresponding to a number $n$ as $fun(n)$. We proved before that for prime $n$ all the Fibonacci numbers that are divisible by $n$ have indices that are multpiples of $fun(n)$. This fact is also true for any $n$ (see Wall \cite{wall} for the proof). Suppose two numbers $n$ and $m$ are mutually coprime. Then the Fibonacci numbers that are divisible by $nm$ have indices that are multiples of both $fun(n)$ and $fun(m)$. The last statement is equivalent to 
$$fun(nm) = \gcd(fun(n),fun(m)).$$
If no odd Fibonacci number is divisible by $nm$, then $fun(nm)$ is divisible by 3. Hence, $fun(n)$ or $fun(m)$ is divisible by 3. That means, if $nm$ is in the complete sequence A133246 for coprime $m$ and $n$, then $m$ or $n$ are in the sequence. Hence $nm$ doesn't belong to the base sequence.

Now it remains to investigate the powers of prime numbers.

\begin{lemma}
For a prime $p > 2$ and a positive integer $t$, $fun(p^{t+1})$ equals either $fun(p^t)$ or $p*fun(p^t)$.
\end{lemma}

\begin{proof}
To prove the lemma I will use the explicit formula for the Fibonacci numbers. Let us denote the golden ratio by $\phi$, then $F_n = (\phi^n - \phi^{-n})/\sqrt{5}$. We will also use Lucas numbers $L_n = \phi^n + \phi^{-n}$ as a helping tool in our calculations. You can check that $L_n = F_{n-1} + F_{n+1}$. The useful thing about $L_n$ is that $\gcd{(F_n,L_n)}$ is either 1 or 2. 

Suppose $m$ is the fundamental period for $p^t$. We know that the fundamental period of $p^{t+1}$ is a multiple of $m$. If $fun(p^{t+1}) = m$, we are done. Otherwise, we would like to look at Fibonacci numbers of the form $F_{am}$, where $a > 1$. Let us express $F_{am}$ in terms of $F_m$ and $L_m$ (here I follow a similar calculation in Wall's paper \cite{wall}):

$$F_{am} = (\phi^{am} - \phi^{-am})/\sqrt{5} = 2^{-a}((\sqrt{5}F_m + L_m)^a - (-\sqrt{5}F_m + L_m)^a)/\sqrt{5}.$$

I am interested in the divisibility of $F_{am}$ by $p^{t+1}$. That means I can ignore the coefficient $2^{-a}$. Also, the non-trivial powers of $F_m$ are divisible by $p^{t+1}$. Hence, I do not need to write the whole expansion formula for $F_{am}$. I am only interested in $((\sqrt{5}F_m + L_m)^a - (-\sqrt{5}F_m + L_m)^a)/\sqrt{5} \pmod{F_m^2}$, which is $2aF_mL_m^{a-1}$. This means, that if $F_m$ is not divisible by $p^{t+1}$, then the smallest number divisible by $p^{t+1}$ is $F_{pm}$.
\end{proof}

\begin{corollary}
For a prime $p > 2$, the $fun(p^t)/fun(p)$ is a power of $p$.
\end{corollary}

Hence, if $p$ is not 3, then 3 divides $fun(p^t)$ iff 3 divides $fun(p)$. In other words, no odd Fibonacci divides $p^t$ iff no odd Fibonacci divides $p$. Hence, for $p \ne 3$, non-trivial powers of $p$ can not belong to the base sequence.

At the conclusion of the proof of the theorem we can remember that there is an odd Fibonacci that divides 3 and no odd Fibonacci divides 9.
\end{proof}

\section{Acknowledgements}

I would like to thank Neil Sloane for maintaining the OEIS (Online Encyclopedia of Integer Sequences) --- an extremely helpful resource in studying integer sequences. I decided to write this paper because the sequences I calculated were not present in the OEIS. I am thankful to my sons, Alexey Radul and Sergei Bernstein, for suggesting that I look at the base sequence. I am thankful to Jane Sherwin and Sue Katz for helping me with English for this paper.


\begin{thebibliography}{9}
\bibitem{Honsberger} Ross Honsberger, Mathematical Gems III, 1997.

\bibitem{MWFN} MathWorld article on Fibonacci Numbers. \emph{http://mathworld.wolfram.com/FibonacciNumber.html}

\bibitem{OEIS} N. J. A. Sloane, Online Encyclopedia of Integer Sequences (OEIS). \emph{http://www.research.att.com/$\sim$njas/sequences/}

\bibitem{wall}
D.D. Wall, Fibonacci Series Modulo $m$, The American Mathematica Monthly, Vol. 67, No. 6, (1960), 525-532. 

\bibitem{Pisano} Wiki article on Pisano period. \emph{http://en.wikipedia.org/wiki/Pisano\_period}

\end{thebibliography}
\end{document}